\def\titlep{Pentagon equation arising 
from state equations of a C$^{*}$-bialgebra}
\documentclass[11pt]{article}
\usepackage{graphicx,ifthen}
\usepackage{amssymb}

\font\germ=eufm10 at12pt

\def\goth#1{\hbox{\germ#1}}

\newcommand{\qed}{\hbox{\rule[-2pt]{3pt}{6pt}}}
\newcommand{\qedh}{\hfill\qed \\}

\newcommand{\vep}{\varepsilon}

\newtheorem{Thm}{Theorem}[section]

\newtheorem{rem}[Thm]{Remark}

\newtheorem{defi}[Thm]{Definition}
\newtheorem{lem}[Thm]{Lemma}

\newtheorem{prop}[Thm]{Proposition}
\newtheorem{prob}[Thm]{Problem}

\newtheorem{fig}[Thm]{Figure}

\def\cal#1{\mathcal #1}
\def\con{{\cal O}_{n}}

\def\edot{=1,\ldots,n}
\def\pr{{\it Proof.}\quad}

\def\co#1{{\cal O}_{#1}}
\def\disp#1{{\displaystyle #1}}
\def\brl{branching law}
\def\bfsnl{{\rm BFS}_{N}(\Lambda)}

\addtocounter{footnote}{1}
\def\sftt#1{
\setcounter{equation}{0}
\addtocounter{footnote}{1}
\section{#1}
}

\def\ssft#1{\subsection{#1}}

\def\ssfr#1{\subsection*{#1}}

\begin{document}
%
%
\def\autherp{Katsunori Kawamura}
\def\emailp{e-mail: kawamura@kurims.kyoto-u.ac.jp.}
\def\addressp{{\small {\it College of Science and Engineering Ritsumeikan University,}}\\
{\small {\it 1-1-1 Noji Higashi, Kusatsu, Shiga 525-8577, Japan}}
}

\def\infw{\Lambda^{\frac{\infty}{2}}V}
\def\zhalfs{{\bf Z}+\frac{1}{2}}
\def\ems{\emptyset}
\def\pmvac{|{\rm vac}\!\!>\!\! _{\pm}}
\def\vac{|{\rm vac}\rangle _{+}}
\def\dvac{|{\rm vac}\rangle _{-}}
\def\ovac{|0\rangle}
\def\tovac{|\tilde{0}\rangle}
\def\expt#1{\langle #1\rangle}
\def\zph{{\bf Z}_{+/2}}
\def\zmh{{\bf Z}_{-/2}}
\def\brl{branching law}
\def\bfsnl{{\rm BFS}_{N}(\Lambda)}
\def\scm#1{S({\bf C}^{N})^{\otimes #1}}
\def\mqb{\{(M_{i},q_{i},B_{i})\}_{i=1}^{N}}
\def\zhalf{\mbox{${\bf Z}+\frac{1}{2}$}}
\def\zmha{\mbox{${\bf Z}_{\leq 0}-\frac{1}{2}$}}
\newcommand{\mline}{\noindent
\thicklines
\setlength{\unitlength}{.1mm}
\begin{picture}(1000,5)
\put(0,0){\line(1,0){1250}}
\end{picture}
\par
 }
\def\ptimes{\otimes_{\varphi}}
\def\qtimes{\otimes_{\tilde{\varphi}}}
\def\delp{\Delta_{\varphi}}
\def\delps{\Delta_{\varphi^{*}}}
\def\gamp{\Gamma_{\varphi}}
\def\gamps{\Gamma_{\varphi^{*}}}
\def\sem{{\sf M}}
\def\sen{{\sf N}}
\def\hdelp{\hat{\Delta}_{\varphi}}
\def\tilco#1{\tilde{\co{#1}}}
\def\ndm#1{{\bf M}_{#1}(\{0,1\})}
\def\fs{{\cal F}{\cal S}({\bf N})}
\def\ba{\mbox{\boldmath$a$}}
\def\bb{\mbox{\boldmath$b$}}
\def\bc{\mbox{\boldmath$c$}}
\def\be{\mbox{\boldmath$e$}}
\def\bp{\mbox{\boldmath$p$}}
\def\bq{\mbox{\boldmath$q$}}
\def\bu{\mbox{\boldmath$u$}}
\def\bv{\mbox{\boldmath$v$}}
\def\bw{\mbox{\boldmath$w$}}
\def\bx{\mbox{\boldmath$x$}}
\def\by{\mbox{\boldmath$y$}}
\def\bz{\mbox{\boldmath$z$}}

%
%
%
\setcounter{section}{0}
\setcounter{footnote}{0}
\setcounter{page}{1}
\pagestyle{plain}

%
%
\title{\titlep}
\author{\autherp\thanks{\emailp}
\\
\addressp}
\date{}
\maketitle
%
%
\begin{abstract}
The direct sum ${\cal O}_{*}$ of all Cuntz algebras
has a non-cocommutative comultiplication $\Delta_{\varphi}$
such that there exists no antipode of any dense subbialgebra 
of the C$^{*}$-bialgebra $({\cal O}_{*},\Delta_{\varphi})$. 
From states equations of ${\cal O}_{*}$
with respect to the tensor product,
we construct an operator $W$ for 
$({\cal O}_{*},\Delta_{\varphi})$
such that  $W^{*}$ is an isometry,
$W(x\otimes I)W^{*}=\Delta_{\varphi}(x)$
for each $x\in {\cal O}_{*}$
and $W$ satisfies the pentagon equation.
\end{abstract}


%
%
\sftt{Introduction}
\label{section:first}
Quantum groups have roots in solvable lattice model 
as mathematical physics \cite{Drinfeld, Jimbo}.
On the other hand, similar objects were studied in operator algebra
as a generalization of the Pontryagin
duality for abelian locally compact groups 
by using C$^{*}$-bialgebras \cite{KV, MNW}.
We have studied C$^{*}$-bialgebras and their representations.
In this paper, 
we construct a kind of multiplicative isometry for
a C$^{*}$-bialgebra from states which satisfy 
tensor product equations induced by the comultiplication.
In this section, we show our motivation,
definitions of C$^{*}$-bialgebras and our main theorem.

%
%
\ssft{Motivation}
\label{subsection:firstone}
In this subsection, we roughly explain our motivation 
and the background of this study.
Explicit mathematical definitions will 
be shown after $\S$ \ref{subsection:firsttwo}.

Define the C$^{*}$-algebra $\co{*}$
as the direct sum of all Cuntz algebras except $\co{\infty}$:
%
%
\begin{equation}
\label{eqn:cuntztwo}
\co{*}=\co{1}\oplus\co{2}\oplus\co{3}\oplus\co{4}\oplus\cdots,
\end{equation}
where $\co{1}$ denotes the $1$-dimensional C$^{*}$-algebra ${\bf C}$
for convenience.
In \cite{TS02}, 
we showed that ${\cal O}_{*}$ has 
a non-cocommutative comultiplication $\Delta_{\varphi}$
such that there exists no antipode of any dense subbialgebra 
of the C$^{*}$-bialgebra $({\cal O}_{*},\Delta_{\varphi})$. 
We investigated a Haar state, KMS states, C$^{*}$-bialgebra automorphisms 
and C$^{*}$-subbialgebras. 
This study was motivated by a certain tensor product 
of representations of Cuntz algebras \cite{TS01}.
With respect to the tensor product, 
tensor product formulae for irreducible representations 
and type III factor representations were computed \cite{TS01,TS07}.
Since there is no standard comultiplication of Cuntz algebras,
$\co{*}$ is not a deformation of any known cocommutative bialgebra. 
The C$^{*}$-bialgebra $\co{*}$ is a rare example 
of not only C$^{*}$-bialgebra but also purely algebraic bialgebra.
Hence we are interested in the bialgebra structure of $\co{*}$. 

On the other hand,
C$^{*}$-bialgebras have been studied in quantum groups in operator algebras
\cite{KV,MNW}.
In order to investigate the C$^{*}$-bialgebra $\co{*}$, 
the theory of quantum groups is one of leading cases even if 
the original motivation of the study of $\co{*}$ is not a quantum group.
Hence we are interested whether various statements of quantum groups 
hold on $\co{*}$ or not.

For example, 
the Kac-Takesaki operator is important to describe the dual of a quantum group.
As a study of duality for groups,
it was introduced by Stinespring \cite{Stinespring}, Kac \cite{Katz01,Katz02} 
and Takesaki \cite{Takesaki},
and was generalized to locally compact quantum groups
by \cite{KV,MNW}.
Furthermore, the Kac-Takesaki operator was generalized to
multiplicative unitary \cite{BS}.
In \cite{KV,MNW}, a C$^{*}$-bialgebra $A$ 
with an invariant weight $\omega$ is considered,
and a antipode and a Kac-Takesaki operator are naturally induced 
from this setting $(A,\omega)$ (\cite{MNW}, Theorem 1.9).
Since $(\co{*},\delp)$ never has antipode, there exists 
no such weight on $\co{*}$.
Hence, in this paper, we consider the following question instead 
of the existence of invariant weight:
%
%
\begin{prob}
\label{prob:first}
Find an operator $W$ for $(\co{*},\delp)$
such that 
\begin{enumerate}
\item
$W(x\otimes I)W^{*}=\delp(x)$ for each $x\in \co{*}$, and 
\item
$W$ satisfies the pentagon equation.
\end{enumerate}
\end{prob}

%
%
\ssft{Covariant representation of C$^{*}$-bialgebra}
\label{subsection:firsttwo}
In this subsection, we recall definitions of C$^{*}$-bialgebra,
and we introduce covariant representation of a C$^{*}$-bialgebra.

At first,
we prepare terminologies about C$^{*}$-bialgebra.
Assume that every tensor product $\otimes$ as below means 
the minimal C$^{*}$-tensor product.
Let $M(A)$ denote the multiplier algebra of $A$.
A pair $(A,\Delta)$ is a {\it C$^{*}$-bialgebra}
if $A$ is a C$^{*}$-algebra and $\Delta$ is a $*$-homomorphism 
from $A$ to $M(A\otimes A)$ 
such that the linear span of $\{\Delta(a)(b\otimes c):a,b,c\in A\}$ 
is norm dense in $A\otimes A$, and the following holds:
%
%
\begin{equation}
\label{eqn:bialgebratwo}
(\Delta\otimes id)\circ \Delta=(id\otimes\Delta)\circ \Delta.
\end{equation}
We call $\Delta$ the {\it comultiplication} of $A$.
We state that a C$^{*}$-bialgebra $(A,\Delta)$ is {\it strictly proper}
if $\Delta(a)\in A\otimes A$ for any $a\in A$.
For two strictly proper  
C$^{*}$-bialgebras $(A_{1},\Delta_{1})$ and 
$(A_{2},\Delta_{2})$,
$f$ is a {\it strictly proper C$^{*}$-bialgebra morphism} from 
$(A_{1},\Delta_{1})$ to $(A_{2},\Delta_{2})$ 
if $f$ is $*$-homomorphism from $A_{1}$ to $A_{2}$
such that $(f\otimes f)\circ \Delta_{1}=\Delta_{2}\circ f$.
In addition, if $f$ is bijective,
then $f$ is called a {\it strictly proper C$^{*}$-bialgebra isomorphism}.
Remark that a locally compact quantum group as a C$^{*}$-bialgebra 
is not always strictly proper \cite{KV,MNW}.

Let $(A,\Delta)$ be a strictly proper C$^{*}$-bialgebra.
If $({\cal H},\pi)$ is 
a faithful $*$-representation of $A$,
then we can define the comultiplication $\Delta^{'}$ on $\pi(A)$ as follows:
%
%
\begin{equation}
\label{eqn:repdelta}
\Delta^{'}\equiv (\pi\otimes \pi)\circ \Delta\circ \pi^{-1}.
\end{equation}
Then $(\pi(A),\Delta^{'})$ is also a 
strictly proper C$^{*}$-bialgebra which is isomorphic to $(A,\Delta)$.

We reformulate Problem \ref{prob:first}
by introducing a representation of a C$^{*}$-bialgebra as follows. 
%
%
\begin{defi}
\label{defi:kto}
Let $(A,\Delta)$ be a strictly proper C$^{*}$-bialgebra.
\begin{enumerate}
\item
A triplet $({\cal H},\pi,W)$ is a quasi-covariant representation of $(A,\Delta)$
if $({\cal H},\pi)$ is a $*$-representation of the C$^{*}$-algebra $A$
and $W$ is a nonzero partial isometry on ${\cal H}\otimes {\cal H}$ such that 
%
%
\begin{equation}
\label{eqn:ktoone}
W (\pi(a) \otimes I)=\{(\pi\otimes \pi)\circ \Delta\}(a)W \quad(a\in A)
\end{equation}
where $I$ denotes the identity operator on ${\cal H}$.
\item
In addition to (i),
if $W^{*}$ is an isometry, then 
we call $({\cal H},\pi,W)$  a covariant representation of $(A,\Delta)$.
\item
A quasi-covariant representation $({\cal H},\pi,W)$ of $(A,\Delta)$ is 
pentagonal if $W$ satisfies the following pentagon equation on the three fold
tensor ${\cal H}\otimes {\cal H}\otimes {\cal H}$ of ${\cal H}$:
%
%
\begin{equation}
\label{eqn:pen}
W_{12}W_{13}W_{23}=W_{23}W_{12}
\end{equation}
where we use the leg numbering notation in \cite{BS}.
\end{enumerate}
\end{defi}

\noindent
If a unitary $W$ satisfies (\ref{eqn:pen}),
then $W$ is called a {\it multiplicative unitary} \cite{BS}.
As generalizations,
a multiplicative isometry and a multiplicative partial isometry
are considered in \cite{BS3}.
%
%
\begin{rem}
\label{rem:second}
{\rm
In Definition \ref{defi:kto},
the choice of $W$ has the ambiguity of the $U(1)$-freedom at least.
From this, 
a quasi-covariant representation is not always pentagonal.
We assume {\it neither} the unitarity {\it nor} 
the pentagon equation for $W$ in general.
Furthermore, we do not consider the uniqueness of a covariant representation
$(A,\Delta)$.
If we identify $\pi(x)$ with $x$ and $W^{*}$ is an isometry,
then (\ref{eqn:ktoone}) is rewritten as follows:
%
%
\begin{equation}
\label{eqn:ktotwo}
\Delta(a) = W (a \otimes I)W^{*}\quad(a\in A).
\end{equation}
}
\end{rem}
For (\ref{eqn:ktotwo}),
it is often said that 
the comultiplication $\Delta$ 
is implemented by $W$ (\cite{MNW}, Proposition 3.6(1)).
In Theorem \ref{Thm:multi}, 
we will construct pentagonal covariant representations 
of $\co{*}$ in (\ref{eqn:cuntztwo}).

%
%
\ssft{C$^{*}$-bialgebra $(\co{*},\delp)$}
\label{subsection:firstthree}
In this subsection, we recall the C$^{*}$-bialgebra in \cite{TS02}.
Let $\con$ denote the Cuntz algebra for $2\leq n<\infty$ \cite{C},
that is, the C$^{*}$-algebra which is universally generated by
generators $s_{1},\ldots,s_{n}$ satisfying
$s_{i}^{*}s_{j}=\delta_{ij}I$ for $i,j\edot$ and
$\sum_{i=1}^{n}s_{i}s_{i}^{*}=I$
where $I$ denotes the unit of $\con$.
The Cuntz algebra $\con$ is simple, that is,
there is no nontrivial two-sided closed ideal.
This implies that any unital representation of $\con$ is faithful.

Redefine the C$^{*}$-algebra $\co{*}$
as the direct sum of the set $\{\con:n\in {\bf N}\}$ of Cuntz algebras:
%
%
\begin{equation}
\label{eqn:cuntbi}
\co{*}\equiv \bigoplus_{n\in {\bf N}} \con
=\{(x_{n}):\|(x_{n})\|\to 0\mbox{ as }n\to\infty\}
\end{equation}
where ${\bf N}=\{1,2,3,\ldots\}$
and $\co{1}$ denotes the $1$-dimensional C$^{*}$-algebra for convenience.
For $n\in {\bf N}$,
let $I_{n}$ denote the unit of $\con$ 
and let $s_{1}^{(n)},\ldots,s_{n}^{(n)}$ denote
canonical generators of $\con$
where $s_{1}^{(1)}\equiv I_{1}$.
For $n,m\in {\bf N}$,
define the embedding $\varphi_{n,m}$ of $\co{nm}$
into $\con\otimes \co{m}$ by
%
%
\begin{equation}
\label{eqn:embeddingone}
\varphi_{n,m}(s_{m(i-1)+j}^{(nm)})\equiv s_{i}^{(n)}\otimes s_{j}^{(m)}
\quad(i=1,\ldots,n,\,j=1,\ldots,m).
\end{equation}
%
%
\begin{Thm}
\label{Thm:mainone}
For the set $\varphi\equiv \{\varphi_{n,m}:n,m\in {\bf N}\}$ in
(\ref{eqn:embeddingone}),
define the $*$-homomorphism $\delp$ from $\co{*}$ to $\co{*}\otimes \co{*}$ by
%
%
\begin{eqnarray}
\label{eqn:dpone}
\delp\equiv& \oplus\{\delp^{(n)}:n\in {\bf N}\},\\
\nonumber
\\
\label{eqn:dptwo}
\delp^{(n)}(x)\equiv &\disp{\sum_{(m,l)\in {\bf N}^{2},\,ml=n}\varphi_{m,l}(x)}
\quad(x\in \con,\,n\in {\bf N}).
\end{eqnarray}
Then  the following holds:
\begin{enumerate}
\item
(\cite{TS02}, Theorem 1.1)
The pair $(\co{*},\delp)$ is a  strictly proper non-cocommutative C$^{*}$-bialgebra.
\item
(\cite{TS02}, Theorem 1.2(v))
 There is no antipode for any dense subbialgebra of $\co{*}$.
\end{enumerate}
\end{Thm}
About properties of $\co{*}$, see \cite{TS02,TS04}.
About a generalization of $\co{*}$, see \cite{TS05}.

Let ${\rm Rep}\con$ denote the class of all $*$-representations of $\con$.
For $\pi_{1},\pi_{2}\in{\rm Rep}\con$,
we define the relation 
$\pi_{1}\sim \pi_{2}$ if $\pi_{1}$ and $\pi_{2}$ are unitarily equivalent.
Then the following holds.
%
%
\begin{lem}(\cite{TS01}, Lemma 1.2)
\label{lem:fundamental}
For $\varphi_{n,m}$ in (\ref{eqn:embeddingone}),
$\pi_{1}\in {\rm Rep}\con$ and  $\pi_{2}\in {\rm Rep}\co{m}$,
define $\pi_{1}\ptimes \pi_{2}\in {\rm Rep}\co{nm}$ by
%
%
\begin{equation}
\label{eqn:ptimes}
\pi_{1}\ptimes \pi_{2}\equiv (\pi_{1}\otimes \pi_{2})\circ \varphi_{n,m}.
\end{equation}
Then the following holds for
$\pi_{1},\pi_{1}^{'}\in {\rm Rep}\con$,
$\pi_{2},\pi_{2}^{'}\in {\rm Rep}\co{m}$ and $\pi_{3}\in {\rm Rep}\co{l}$:
\begin{enumerate}
\item
If $\pi_{1}\sim \pi_{1}^{'}$ and $\pi_{2}\sim \pi_{2}^{'}$,
then $\pi_{1}\otimes_{\varphi} \pi_{2}\sim
\pi_{1}^{'}\otimes_{\varphi} \pi_{2}^{'}$.
\item
$\pi_{1}\otimes_{\varphi} (\pi_{2}\oplus \pi_{2}^{'})=
\pi_{1}\otimes_{\varphi} \pi_{2}\,\oplus\, \pi_{1}\otimes_{\varphi} \pi_{2}^{'}$.
\item
$\pi_{1}\otimes_{\varphi} (\pi_{2}\otimes_{\varphi} \pi_{3})
=(\pi_{1}\otimes_{\varphi} \pi_{2})\otimes_{\varphi} \pi_{3}$.
\end{enumerate}
\end{lem}

\noindent
From Lemma \ref{lem:fundamental}(i),
we can define $[\pi_{1}]\ptimes [\pi_{2}]\equiv [\pi_{1}\ptimes \pi_{2}]$
where  $[\pi]$  denotes the unitary equivalence class of $\pi$.

Let ${\cal S}_{n}$ denote the set of all states of $\con$.
For $(\omega,\omega^{'})\in {\cal S}_{n}\times {\cal S}_{m}$,
define
%
%
\begin{equation}
\label{eqn:eleven}
\omega\ptimes \omega^{'}\equiv (\omega\otimes \omega^{'})\circ
\varphi_{n,m}
\end{equation}
where $(\omega\otimes \omega^{'})(x\otimes y)\equiv \omega(x)\omega^{'}(y)$
for $x\in \con$ and $y\in\co{m}$.
Then we see that 
$\omega \ptimes (\omega^{'}\ptimes \omega^{''})
=
(\omega \ptimes \omega^{'})\ptimes \omega^{''}$.

%
%
\ssft{Main theorem}
\label{subsection:firstfour}
In this subsection, we show our main theorem.
%
%
\begin{Thm}
\label{Thm:multi}
Assume that $\{\omega_{n}:n\geq 1\}$ is a set of states
such that $\omega_{n}$ is a state of $\con$ 
with the Gel'fand-Na\u{\i}mark-Segal (=GNS) 
triple $({\cal H}_{n},\pi_{n},\Omega_{n})$ for $n\geq 1$ and 
%
%
\begin{equation}
\label{eqn:statemultib}
\omega_{n}\ptimes \omega_{m}=\omega_{nm}\quad(n,m\geq 1)
\end{equation}
where $\ptimes$ is as in (\ref{eqn:eleven}).
Then there exists
a nonzero partial isometry $W^{(n,m)}$ from ${\cal H}_{nm}\otimes {\cal H}_{m}$
to ${\cal H}_{n}\otimes {\cal H}_{m}$ for each $n,m\geq 1$ 
such that the following holds:
\begin{enumerate}
\item 
For each $n,m\geq 1$,
%
%
\begin{equation}
\label{eqn:wzy}
W^{(n,m)}(\pi_{nm}(X)\otimes I_{m})=
(\pi_{n}\ptimes \pi_{m})(X)W^{(n,m)}\quad(X\in\co{nm})
\end{equation}
where $I_{m}$ denotes the identity operator on ${\cal H}_{m}$
and $\ptimes$ is as in (\ref{eqn:ptimes}).
\item
In addition,
if\, $\Omega_{n}\otimes \Omega_{m}$ is a cyclic vector for 
$({\cal H}_{n}\otimes {\cal H}_{m},\pi_{n}\ptimes \pi_{m})$, 
then we can choose $W^{(n,m)}$ such that  $(W^{(n,m)})^{*}$ is an isometry.
\item
For each $n,m,l\geq 1$,
%
%
\begin{equation}
\label{eqn:urel}
W^{(n,m)}_{12}W^{(nm,l)}_{13}W^{(m,l)}_{23}=
W^{(m,l)}_{23}W^{(n,ml)}_{12}
\end{equation}
on ${\cal H}_{nml}
\otimes {\cal H}_{ml}\otimes {\cal H}_{l}$.
\item
Let $(\co{*},\delp)$ be as in Theorem \ref{Thm:mainone}.
Define
%
%
\begin{equation}
\label{eqn:main}
{\cal H}\equiv \bigoplus_{n\in{\bf N}}{\cal H}_{n},\quad
\pi\equiv \bigoplus_{n\in{\bf N}}\pi_{n},\quad
W\equiv \bigoplus_{n,m\in {\bf N}}W^{(n,m)}.
\end{equation}
Then $({\cal H},\pi,W)$ is a  
pentagonal quasi-covariant representation of $(\co{*},\delp)$.
In addition,
if the assumption in (ii) is satisfied for each $n,m\geq 1$,
then 
$({\cal H},\pi,W)$ is a pentagonal covariant representation of $(\co{*},\delp)$.
\end{enumerate}
\end{Thm}

\noindent
In consequence,
we obtain a solution of Problem \ref{prob:first}
when a set of states in (\ref{eqn:statemultib}) is given.
The equation (\ref{eqn:urel}) will be generalized 
and closely explained in $\S$ \ref{subsection:secondtwo}.
%
%
\begin{rem}
\label{rem:one}
{\rm
\begin{enumerate}
\item
The operator $W$ in (1.16) does not satisfy the axiom
of multiplicative partial isometry in $\S$ 2 of \cite{BS3}.
\item
The assumption in (ii) does not always hold even if 
(\ref{eqn:statemultib}) holds.
\item
A relation between Cuntz algebras and multiplicative unitaries
is studied by Roberts \cite{Roberts},
which is different from our use.
\end{enumerate}
}
\end{rem}

%
%
\begin{prob}
\label{prob:one}
{\rm
\begin{enumerate}
\item
Generalize Theorem \ref{Thm:multi} to general C$^{*}$-bialgebras.
\item
Show a duality type theorem for $(\co{*},\delp)$.
\item
Dose there exist the pentagonal covariant representation $({\cal H},\pi,W)$
of $(\co{*},\delp)$ such that $W$ is a unitary?
\end{enumerate}
}
\end{prob}

In $\S$ \ref{section:second},
we will show general results and
prove Theorem \ref{Thm:multi}.
In $\S$ \ref{section:third},
we will show examples of states which satisfy equations 
in (\ref{eqn:statemultib})
and the assumption in Theorem \ref{Thm:multi}(ii).

%
%
\sftt{Proof of Theorem \ref{Thm:multi}}
\label{section:second}
In order to prove Theorem \ref{Thm:multi}, 
we show general statements about C$^{*}$-bialgebras in this section.	
%
%
\ssft{C$^{*}$-weakly coassociative system}
\label{subsection:secondone}
We review C$^{*}$-weakly coassociative system in $\S$ 3 of \cite{TS02}.
A {\it monoid} is a set $\sem$ equipped with a binary associative operation 
$\sem\times \sem\ni(a,b)\mapsto ab\in \sem$
and a unit with respect to the operation.
%
%
\begin{defi}
\label{defi:axiom}
Let $\sem$ be a monoid with a unit $e$.
A data $(\{A_{a}:a\in \sem\},\{\varphi_{a,b}:a,b\in \sem\})$
is a C$^{*}$-weakly coassociative system (= C$^{*}$-WCS) over $\sem$ if 
$A_{a}$ is a unital C$^{*}$-algebra with a unit $I_{a}$ for $a\in \sem$
and $\varphi_{a,b}$ is a unital $*$-homomorphism
from $A_{ab}$ to $A_{a}\otimes A_{b}$
for $a,b\in \sem$ such that
\begin{enumerate}
\item
for all $a,b,c\in \sem$, the following holds:
%
%
\begin{equation}
\label{eqn:wcs}
(id_{a}\otimes \varphi_{b,c})\circ \varphi_{a,bc}
=(\varphi_{a,b}\otimes id_{c})\circ \varphi_{ab,c}
\end{equation}
where $id_{x}$ denotes the identity map on $A_{x}$ for $x=a,c$,
\item
there exists a counit $\vep_{e}$ of $A_{e}$ 
such that $(A_{e},\varphi_{e,e},\vep_{e})$ is a counital C$^{*}$-bialgebra,
\item
$\varphi_{e,a}(x)=I_{e}\otimes x$ and
$\varphi_{a,e}(x)=x\otimes I_{e}$ for $x\in A_{a}$ and $a\in \sem$.
\end{enumerate}
\end{defi}

\noindent
The system $(\{\con:n\in {\bf N}\},\{\varphi_{n,m}:n,m\in {\bf N}\})$
in (\ref{eqn:embeddingone}) is a C$^{*}$-WCS.
As for the other example of C$^{*}$-WCS, see $\S$ 1.3 of \cite{TS05}.
%
%
\begin{Thm}(\cite{TS02}, Theorem 3.1)
\label{Thm:subthree}
Let $(\{A_{a}:a\in \sem\},\{\varphi_{a,b}:a,b\in \sem\})$ be a C$^{*}$-WCS 
over a monoid $\sem$.
Assume that $\sem$ satisfies that 
%
%
\begin{equation}
\label{eqn:finiteness}
\#{\cal N}_{a}<\infty \mbox{ for each }a\in \sem
\end{equation}
where ${\cal N}_{a}\equiv\{(b,c)\in \sem\times \sem:\,bc=a\}$.
Define the C$^{*}$-algebra 
%
%
\begin{equation}
\label{eqn:astar}
A_{*}\equiv  \oplus \{A_{a}:a\in \sem\},
\end{equation}
and define the $*$-homomorphism $\Delta_{\varphi}$ from $A_{*}$ 
to $A_{*}\otimes A_{*}$ by
%
%
\begin{equation}
\label{eqn:directsumtwo}
\Delta_{\varphi}\equiv \oplus\{\Delta_{\varphi}^{(a)}:a\in \sem\},
\quad \Delta^{(a)}_{\varphi}(x)\equiv \sum_{(b,c)\in {\cal N}_{a}}
\varphi_{b,c}(x)\quad(x\in A_{a}).
\end{equation}
Then $(A_{*},\delp)$ is a strictly proper C$^{*}$-bialgebra.
\end{Thm}

\noindent
We call $(A_{*},\Delta_{\varphi})$ the {\it C$^{*}$-bialgebra associated with} 
$(\{A_{a}:a\in \sem\},\{\varphi_{a,b}:a,b\in \sem\})$.
In this paper,
we always assume (\ref{eqn:finiteness}).

Let ${\rm Rep}A_{a}$ denote the class of all $*$-representations of $A_{a}$.
For $\pi_{a}\in {\rm Rep}A_{a}$
and $\pi_{b}\in {\rm Rep}A_{b}$,
define
$\pi_{a}\ptimes \pi_{b}\in {\rm Rep}A_{ab}$
by
%
%
\begin{equation}
\label{eqn:ptensor}
\pi_{a}\ptimes \pi_{b}\equiv 
(\pi_{a}\otimes \pi_{b})\circ \varphi_{a,b}.
\end{equation}
From (\ref{eqn:wcs}),
we see that statements in Lemma \ref{lem:fundamental} also hold for $\ptimes$
in (\ref{eqn:ptensor}).

%
%
\ssft{Covariant representation of C$^{*}$-WCS}
\label{subsection:secondtwo}
We introduce covariant representation of C$^{*}$-WCS in this subsection.
%
%
\begin{defi}
\label{defi:cov}
\begin{enumerate}
\item
A data $(\{({\cal H}_{a},\pi_{a}):a\in\sem\},\,\{W^{(a,b)}:a,b\in\sem\})$
is a quasi-covariant representation of a C$^{*}$-WCS 
$(\{A_{a}:a\in \sem\},\{\varphi_{a,b}:a,b\in\sem\})$ 
if $({\cal H}_{a},\pi_{a})$ is a unital $*$-representation 
of the C$^{*}$-algebra $A_{a}$ and 
$W^{(a,b)}$ is a nonzero partial isometry from 
${\cal H}_{ab}\otimes {\cal H}_{b}$ to
${\cal H}_{a}\otimes {\cal H}_{b}$
such that $W^{(a,b)}$ satisfies 
%
%
\begin{equation}
\label{eqn:covone}
W^{(a,b)}(\pi_{ab}(x)\otimes I_{b})
=(\pi_{a}\ptimes \pi_{b})(x)W^{(a,b)}\quad(x\in A_{ab})
\end{equation}
for each $a,b\in \sem$
where $I_{b}$ denotes the identity operator on ${\cal H}_{b}$.
\item
In addition to (i),
if $(W^{(a,b)})^{*}$ is an isometry,
we call $(\{({\cal H}_{a},\pi_{a}):a\in\sem\},\,\{W^{(a,b)}:a,b\in\sem\})$
a covariant representation of 
$(\{A_{a}:a\in \sem\},\{\varphi_{a,b}:a,b\in\sem\})$.
\item
A quasi-covariant representation 
$(\{({\cal H}_{a},\pi_{a}):a\in\sem\},\,\{W^{(a,b)}:a,b\in\sem\})$
is pentagonal if the following relation holds
on ${\cal H}_{abc}\otimes {\cal H}_{bc}\otimes {\cal H}_{c}$
for each $a,b,c\in\sem$:
%
%
\begin{equation}
\label{eqn:cwcsp}
W^{(a,b)}_{12}W^{(ab,c)}_{13}W^{(b,c)}_{23}
=W^{(b,c)}_{23}W^{(a,bc)}_{12}.
\end{equation}
\end{enumerate}
\end{defi}
We illustrate (\ref{eqn:cwcsp}) as the commutative diagram 
in Figure \ref{fig:diagram}:

\def\hilbert#1#2#3{${\cal H}_{#1}\otimes {\cal H}_{#2}\otimes {\cal H}_{#3}$}
\def\figone{
\put(350,450){\hilbert{abc}{bc}{c}}
\put(-50,250){\hilbert{abc}{b}{c}}
\put(-50,0){\hilbert{ab}{b}{c}}
\put(750,250){\hilbert{a}{bc}{c}}
\put(780,00){\hilbert{a}{b}{c}}
\put(300,400){\vector(-2,-1){170}}
\put(650,400){\vector(2,-1){170}}
\put(100,200){\vector(0,-1){130}}
\put(900,200){\vector(0,-1){130}}
\put(300,10){\vector(1,0){400}}
\put(80,370){$W^{(b,c)}_{23}$}
\put(-50,120){$W^{(ab,c)}_{13}$}
\put(450,50){$W^{(a,b)}_{12}$}
\put(790,370){$W^{(a,bc)}_{12}$}
\put(930,120){$W^{(b,c)}_{23}$}
}
%
%
\begin{fig}
\label{fig:diagram}
\quad\\
\thicklines
\setlength{\unitlength}{.1mm}
\begin{picture}(1200,500)(-120,-50)
\put(0,0){\figone}
\end{picture}
\end{fig}

\noindent
where
$W^{(ab,c)}_{13}$ means
$\tau_{2,3}^{-1}\circ (W^{(ab,c)}\otimes I_{b})\circ \tau_{2,3}$
and $\tau_{2,3}$ denotes the permutation of the second Hilbert space 
and the third one.

A monoid $\sem$ is {\it cancellative} 
if the following is satisfied for each $b,c,b^{'}\in\sem$:
If $bc=b^{'}c$, then $b=b^{'}$,
and if $cb=cb^{'}$, then $b=b^{'}$
(\cite{Petrich}, p 6).
%
%
\begin{prop}
\label{prop:cov}
Let $(\{({\cal H}_{a},\pi_{a}):a\in\sem\},\,\{W^{(a,b)}:a,b\in\sem\})$
be a quasi-covariant representation ({\it resp.} a covariant representation) 
of a C$^{*}$-WCS $(\{A_{a}:a\in \sem\},\{\varphi_{a,b}:a,b\in\sem\})$
and assume that $\sem$ is cancellative.
Define
%
%
\begin{equation}
\label{eqn:big}
{\cal H}\equiv \bigoplus_{a\in\sem}{\cal H}_{a},\quad
\pi\equiv \bigoplus_{a\in\sem}\pi_{a},\quad
W\equiv \bigoplus_{a,b\in\sem}W^{(a,b)}.
\end{equation}
Then $({\cal H},\pi,W)$
is a quasi-covariant representation 
({\it resp.} a covariant representation) of $(A_{*},\delp)$.
In addition,
if 
$(\{({\cal H}_{a},\pi_{a}):a\in\sem\},\,\{W^{(a,b)}:a,b\in\sem\})$
is pentagonal, then $({\cal H},\pi,W)$ is also pentagonal.
\end{prop}
%
%
\pr
Since images of $\{W^{(a,b)}\}$ are mutually orthogonal
and $\bigoplus_{a,b}{\cal H}_{a}\otimes 
{\cal H}_{b}={\cal H}\otimes {\cal H}$,
$W$ is a partial isometry. 
Especially, if $(W^{(a,b)})^{*}$ is an isometry for each $a,b$,
then $W^{*}$ is also an isometry. 
Define $W^{(a)}\equiv \bigoplus_{bc=a}W^{(b,c)}$.
From (\ref{eqn:covone}),
we can verify that 
%
%
\begin{equation}
\label{eqn:index}
W^{(a)}(\pi(x)\otimes I)
=(\pi\otimes \pi)(\delp(x))W^{(a)}\quad(x\in A_{a}).
\end{equation}
This implies the first statement.

Assume that (\ref{eqn:cwcsp}) is satisfied.
It is sufficient to show that the pentagon equation of $W$
holds on ${\cal H}_{a}\otimes{\cal H}_{b}\otimes{\cal H}_{c}$
for each $a,b,c\in\sem$.
Let $v\in {\cal H}_{a}\otimes{\cal H}_{b}\otimes{\cal H}_{c}$.
Then 
$W_{12}W_{13}W_{23}v=0$ if not $b=b^{'}c$ and $a=a^{'}b^{'}c$
for some $a^{'},b^{'}\in\sem$.
Hence we can assume that  
$v\in {\cal H}_{abc}\otimes{\cal H}_{bc}\otimes{\cal H}_{c}$.
Then we see that
%
%
\begin{equation}
\label{eqn:www}
W_{12}W_{13}W_{23}v=
W_{12}^{(a,b)}W_{13}^{(ab,c)}W_{23}^{(b,c)}v,\quad
W_{23}W_{12}v=W_{23}^{(b,c)}W_{12}^{(a,bc)}v.
\end{equation}
From (\ref{eqn:cwcsp}),
the second statement holds.
\qedh

\noindent
Remark that $W$ in (\ref{eqn:big}) is not a unitary
even if $W^{(a,b)}$ is unitary for each $a,b$,
because
%
%
\begin{equation}
\label{eqn:ker}
{\rm Ker}W=\bigoplus_{b\nmid a}{\cal H}_{a}\otimes {\cal H}_{b} \ne \{0\}
\end{equation}
where the direct sum is taken over all pairs $(a,b)$
such that $b$ is not a right divisor of $a$ in $\sem$.

%
%
\ssft{Multiplicative partial isometry arising 
from states equations for C$^{*}$-WCS}
\label{subsection:secondthree}
In this subsection,
we show that certain tensor equations of states
induce a covariant representation of C$^{*}$-WCS
and prove Theorem \ref{Thm:multi}.

Let $(\{A_{a}:a\in\sem\},\{\varphi_{a,b}:a,b\in\sem\})$
be a C$^{*}$-WCS.
Let $\omega_{a}$ and $\omega_{b}$ be states of $A_{a}$  and $A_{b}$,
respectively.
Define the new state $\omega_{a}\ptimes \omega_{b}$ of 
$A_{ab}$ by
%
%
\begin{equation}
\label{eqn:thirteen}
\omega_{a}\ptimes \omega_{b}
\equiv (\omega_{a}\otimes \omega_{b})\circ \varphi_{a,b}.
\end{equation}
%
%
\begin{lem}
\label{lem:multi}
Let $(\{A_{a}:a\in\sem\},\{\varphi_{a,b}:a,b\in\sem\})$
be a C$^{*}$-WCS, and 
let $\omega_{a}$ be a state of $A_{a}$ 
with the GNS triple $({\cal H}_{a},\pi_{a},\Omega_{a})$ for $a\in\sem$.
Assume 
%
%
\begin{equation}
\label{eqn:fourteen}
\omega_{a}\ptimes \omega_{b}=\omega_{ab}\quad(a,b\in \sem).
\end{equation}
Define the operator $W^{(a,b)}$ from ${\cal H}_{ab}\otimes {\cal H}_{b}$
to ${\cal H}_{a}\otimes {\cal H}_{b}$ by
%
%
\begin{equation}
\label{eqn:woperator}
W^{(a,b)}(\pi_{ab}(x)\Omega_{ab}\otimes v)
\equiv 
(\pi_{a}\otimes \pi_{b})(\varphi_{a,b}(x))(\Omega_{a}\otimes E_{b}v)
\end{equation}
for $x\in A_{ab}$ and $v\in {\cal H}_{b}$
where $E_{b}$ denotes the projection from ${\cal H}_{b}$ 
onto ${\bf C}\Omega_{b}$.
Then the following holds:
\begin{enumerate}
\item
The data $(\{({\cal H}_{a},\pi_{a}):a\in\sem\},\{W^{(a,b)}:a,b\in\sem\})$
is a pentagonal quasi-covariant representation of 
$(\{A_{a}:a\in\sem\},\{\varphi_{a,b}:a,b\in\sem\})$.
%
\item
If 
$\Omega_{a}\otimes \Omega_{b}$ is a cyclic vector for 
$({\cal H}_{a}\otimes {\cal H}_{b},\pi_{a}\ptimes \pi_{b})$
for each $a,b\in\sem$,
then
$(\{({\cal H}_{a},\pi_{a}):a\in\sem\},\{W^{(a,b)}:a,b\in\sem\})$
is a pentagonal covariant representation of 
$(\{A_{a}:a\in\sem\},\{\varphi_{a,b}:a,b\in\sem\})$.
\end{enumerate}
\end{lem}
%
%
\pr
(i)
Let ${\cal K}$ denote the closure of 
$(\pi_{a}\ptimes \pi_{b})(A_{ab})(\Omega_{a}\otimes\Omega_{b})$
in ${\cal H}_{a}\otimes {\cal H}_{b}$.
Then the subrepresentation $(\pi_{a}\ptimes \pi_{b})|_{{\cal K}}$
is unitarily equivalent to $\pi_{ab}$.
Define the isometry $U$ from ${\cal H}_{ab}$ to 
${\cal H}_{a}\otimes {\cal H}_{b}$ by
%
%
\begin{equation}
\label{eqn:double}
U\pi_{ab}(x)\Omega_{ab}\equiv (\pi_{a}\ptimes \pi_{b})(x)
(\Omega_{a}\otimes \Omega_{b})\quad(x\in A_{ab}).
\end{equation}
Then $U$ is well-defined such that
%
%
\begin{equation}
\label{eqn:triple}
U^{*}(\pi_{a}\ptimes \pi_{b})(x)U=\pi_{ab}(x)\quad(x\in A_{ab}).
\end{equation}
Define another isometry $V$ from ${\cal H}_{ab}$
to ${\cal H}_{ab}\otimes {\cal H}_{b}$ by
%
%
\begin{equation}
\label{eqn:quartet}
Vv\equiv v\otimes \Omega_{b}\quad (v\in {\cal H}_{ab}).
\end{equation}
Then
%
%
\begin{equation}
\label{eqn:ppp}
V^{*}(\pi_{ab}(x)\otimes I_{b})V=\pi_{ab}(x)\quad(x\in A_{ab}).
\end{equation}
Since $W^{(a,b)}=UV^{*}$, 
$W^{(a,b)}$ is well-defined and (\ref{eqn:covone}) is satisfied.

By definition,
it is sufficient to show (\ref{eqn:cwcsp})
on the subspace $\pi_{abc}(A_{abc})\Omega_{abc}\otimes {\bf C}
\Omega_{bc}\otimes {\bf C}\Omega_{c}$
of ${\cal H}_{abc}
\otimes {\cal H}_{ab}\otimes {\cal H}_{c}$.
Define 
$\Omega_{a,b,c}\equiv \Omega_{a}\otimes \Omega_{b}\otimes \Omega_{c}$
for $a,b,c\in \sem$.
For $x\in A_{abc}$,\\
$W^{(a,b)}_{12}
W^{(ab,c)}_{13}
W^{(b,c)}_{23}(\pi_{abc}(x)\otimes I_{bc}\otimes I_{c})\Omega_{abc,bc,c}$
%
%
\begin{equation}
\label{eqn:gtwo}
\begin{array}{rl}
=&
W^{(a,b)}_{12}
W^{(ab,c)}_{13}
(\pi_{abc}(x)\otimes I_{b}\otimes I_{c})\Omega_{abc,b,c}\\
=&
W^{(a,b)}_{12}
(\pi_{ab}\otimes \pi_{b}\otimes \pi_{c})((\varphi_{ab,c})_{13}(x))
\Omega_{ab,b,c}\\
=&
\{(\pi_{a}\otimes \pi_{b}\otimes \pi_{c})\circ (\varphi_{a,b}
\otimes id_{c})\circ \varphi_{ab,c}\}(x)\Omega_{a,b,c},\\
\end{array}
\end{equation}

\noindent
$
W^{(b,c)}_{23}
W^{(a,bc)}_{12}(\pi_{abc}(x)\otimes I_{bc}\otimes I_{c})\Omega_{abc,bc,c}$
%
%
\begin{equation}
\label{eqn:gone}
\begin{array}{rl}
=&W^{(b,c)}_{23}
\{(\pi_{a}\otimes \pi_{bc})(\varphi_{a,bc}(x))\otimes I_{c}\}\Omega_{a,bc,c}\\
=&\{(\pi_{a}\otimes \pi_{b}\otimes \pi_{c})\circ 
(id_{a}\otimes \varphi_{b,c})\circ \varphi_{a,bc}\}(x)\Omega_{a,b,c}
\end{array}
\end{equation}
where $(\varphi_{ab,c})_{13}(x)
\equiv (id_{ab}\otimes \tau_{2,3})(\varphi_{ab,c}(x)\otimes I_{b})$
and $\tau_{2,3}$ denotes the permutation of the second component 
and the third one of the tensor product of algebras.
Applying (\ref{eqn:wcs}) to
(\ref{eqn:gtwo}) and (\ref{eqn:gone}),
(\ref{eqn:cwcsp}) holds.

\noindent
(ii)
In the proof of (i),
the operator $U$ is a unitary from the assumption.
Hence $(W^{(a,b)})^{*}=VU^{*}$ is an isometry.
\qedh

\noindent
{\it Proof of Theorem \ref{Thm:multi}.}
Applying Lemma \ref{lem:multi} to 
the C$^{*}$-WCS $(\{\con:n\in{\bf N}\},\{\varphi_{n,m}:n,m\in{\bf N}\})$,
(i), (ii) and (iii) hold. 
Applying Proposition \ref{prop:cov}
to statements in (i), (ii) and (iii), (iv) holds.
\qedh

%
%
\sftt{Pure states of Cuntz algebras parametrized by unit vectors}
\label{section:third}
In this section, we show examples of set of 
states which satisfies (\ref{eqn:statemultib}).
We recall certain states in \cite{GP0123}
and show tensor product formulae among them.
Let $S({\bf C}^{n})$ denote the set $\{z\in {\bf C}^{n}:\|z\|=1\}$
of all unit vectors in ${\bf C}^{n}$.
%
%
\begin{defi}
\label{defi:firstb}(\cite{GP0123}, Proposition 3.1)
For $n\geq 2$,
let $s_{1},\ldots,s_{n}$ denote canonical generators of $\con$.
For $z=(z_{1},\ldots,z_{n})\in S({\bf C}^{n})$,
define the state $\varrho_{z}$ of $\con$ by
%
%
\begin{equation}
\label{eqn:gpstate}
\varrho_{z}(s_{j_{1}}\cdots s_{j_{a}}s_{k_{b}}^{*}\cdots s_{k_{1}}^{*})\equiv 
\overline{z}_{j_{1}}\cdots \overline{z}_{j_{a}}z_{k_{b}}\cdots z_{k_{1}}
\end{equation}
for each $j_{1},\ldots,j_{a},k_{1},\ldots,k_{b}\in \{1,\ldots,n\}$
and $a,b\geq 1$.
\end{defi}
Remark that $a$ and $b$ may not equal in (\ref{eqn:gpstate}).
The following results for $\varrho_{z}$ are known:
For any $z$, $\varrho_{z}$ is pure when $n\geq 2$. 
If $n=1$, we define $\varrho_{z}(x)\equiv x$ for $x\in \co{1}$ 
if and only if $z=1\in S({\bf C}^{1})=U(1)$.
If $z,y\in S({\bf C}^{n})$ and $z\ne y$,
then GNS representations associated with $\varrho_{z}$ and $\varrho_{y}$
are not unitarily equivalent.

Let $\alpha^{(n)}$ denote the canonical $U(n)$-action on $\con$,
that is,
$\alpha_{g}^{(n)}(s_{i})\equiv \sum_{j=1}^{n}g_{ji}s_{j}$
for $i=1,\ldots,n,\,g=(g_{ij})\in U(n)$.
Then the following holds (\cite{TS01}, Proposition 3.1(iii)):  
%
%
\begin{equation}
\label{eqn:agtwo}
(\alpha_{g}^{(n)}\otimes \alpha_{h}^{(m)})\circ \varphi_{n,m}=
\varphi_{n,m}\circ \alpha^{(nm)}_{g\boxtimes h}
\quad(g\in U(n),\,h\in U(m))
\end{equation}
where $g\boxtimes h\in U(nm)$ is defined as
$(g\boxtimes h)_{m(i-1)+j,m(i^{'}-1)+j^{'}}=g_{ii^{'}}h_{jj^{'}}$
for $i,i^{'}=1,\ldots,n$ and $j,j^{'}=1,\ldots,m$.

Let $GP(z)$ denote
the unitary equivalence class of the GNS representation associated 
with $\varrho_{z}$.
For $g\in U(n)$ and the representative $\pi$ of $GP(z)$,
we write $GP(z)\circ \alpha_{g}^{(n)}$ as $[\pi\circ \alpha_{g}^{(n)}]$.
Then the following holds:
%
%
\begin{equation}
\label{eqn:coa}
GP(z)\circ \alpha_{g^{-1}}^{(n)}=GP(gz)\quad(z\in S({\bf C}^{n}),\,g\in U(n))
\end{equation}
where $gz$ denotes the standard action of $U(n)$ on ${\bf C}^{n}$.
Especially,
$GP(1,0,\ldots,0)$ is $P_{n}(1)$ in Definition 1.4(ii) of \cite{TS04}.
If $\pi_{1}$ and $\pi_{2}$
are representatives of $GP(z)$ and $GP(y)$ for 
$z\in S({\bf C}^{n})$ and $y\in S({\bf C}^{m})$, respectively,
then we write $GP(z)\ptimes GP(y)$ as $[\pi_{1}]\ptimes [\pi_{2}]$
for simplicity of description. 
%
%
\begin{Thm}
\label{Thm:commuteb}
For $\ptimes$ in (\ref{eqn:ptimes}),
the following holds
for each $z\in S({\bf C}^{n})$ and $y\in S({\bf C}^{m})$:
\begin{enumerate}
\item
$\varrho_{z}\ptimes \varrho_{y}=\varrho_{z\boxtimes  y}$,
\item
$GP(z)\ptimes GP(y)=GP(z\boxtimes y)$
\end{enumerate}
where $z\boxtimes y\in S({\bf C}^{nm})$ is defined as 
%
%
\begin{equation}
\label{eqn:zy}
(z\boxtimes y)_{m(i-1)+j}\equiv z_{i}y_{j}\quad(i=1,\ldots,n,\,j=1,\ldots,m),
\end{equation}
and we choose $x=1$ when $x\in S({\bf C}^{1})$ for $x=y,z$.
\end{Thm}
%
%
\pr
(i)
By definition, the statement is verified directly. 

\noindent
(ii)
Let $\eta_{n}\equiv (1,0,\ldots,0)\in {\bf C}^{n}\cap S({\bf C}^{n})$.
Then $\eta_{n}\boxtimes \eta_{m}=\eta_{nm}$.
We write $P_{n}(1)\equiv GP(\eta_{n})$.
Choose $g\in U(n)$ and $h\in U(m)$
such that $gz=\eta_{n}$ and $hy=\eta_{m}$.
From these,
$(g^{-1}\boxtimes h^{-1})(\eta_{mn})=z\boxtimes y$.
From (\ref{eqn:coa}),
%
%
\begin{equation}
\label{eqn:gpcov}
GP(z)=P_{n}(1)\circ \alpha_{g}^{(n)},\quad
GP(y)=P_{m}(1)\circ \alpha_{h}^{(m)},\quad 
GP(z\boxtimes y)=P_{nm}(1)\circ \alpha_{g\boxtimes h}^{(nm)}.
\end{equation}
From these and $P_{n}(1)\ptimes P_{m}(1)=P_{nm}(1)$
by Example 4.1 in \cite{TS01},
the statement is verified.
\qedh

%
%
\begin{Thm}
\label{Thm:para}
Assume that a sequence $(z^{(n)})_{n\geq 1}$ satisfies
the following conditions:
%
%
\begin{equation}
\label{eqn:sequence}
z^{(n)}\in S({\bf C}^{n})\quad (n\geq 1),\quad
z^{(n)}\boxtimes z^{(m)}=z^{(nm)}\quad(n,m\geq 1).
\end{equation}
\begin{enumerate}
\item
Then there exists a pentagonal covariant representation
$({\cal H},\pi,W)$ of $(\co{*},\delp)$.
\item
Let $(y^{(n)})_{n\geq 1}$ be another sequence which 
satisfies (\ref{eqn:sequence})
and 
let $({\cal H}^{'},\pi^{'},W^{'})$ be 
the covariant representation corresponding to $(y^{(n)})_{n\geq 1}$.
Then $({\cal H},\pi)$ and 
$({\cal H}^{'},\pi^{'})$ are not unitarily equivalent
when $(z^{(n)})_{n\geq 1}\ne (y^{(n)})_{n\geq 1}$.
\end{enumerate}
\end{Thm}
%
%
\pr
(i)
Define $\omega_{n}\equiv \varrho_{z^{(n)}}$ for $n\geq 1$.
From Theorem \ref{Thm:commuteb}(i) and (\ref{eqn:sequence}),
$\{\omega_{n}:n\geq 1\}$ satisfies 
$\omega_{n}\ptimes \omega_{m}=\omega_{nm}$ for each $n,m\geq 1$.
Since $\omega_{n}$ is pure for each $n$,
the tensor product of GNS representations associated with
$\omega_{n}$ and $\omega_{m}$ is irreducible
from Theorem \ref{Thm:commuteb}(ii).
Therefore the assumption of 
Theorem \ref{Thm:multi}(ii) holds.
From Theorem \ref{Thm:multi}(iv), the statement holds.

\noindent
(ii)
Since $z\in S({\bf C}^{n})$ is a complete invariant of $GP(z)$,
the statement holds from Lemma \ref{lem:multi}.
\qedh

In Theorem \ref{Thm:para},
the crucial point is the choice of sequence $(z^{(n)})_{n\geq 1}$
of unit vectors, which is a monoid with respect to
the product in (\ref{eqn:zy}).
For example,
the following sequences $(z^{(n)})_{n\geq 1}$ satisfy (\ref{eqn:sequence}):
\begin{enumerate}
\item
For $n\geq 1$, $z^{(n)}\equiv (1,0,\ldots,0)\in {\bf C}^{n}$.
\item
For $n\geq 1$, $z^{(n)}\equiv (0,\ldots,0,1)\in {\bf C}^{n}$.
\item
For $n\geq 1$,
$z^{(n)}\equiv (n^{-1/2},\ldots,n^{-1/2})\in {\bf C}^{n}$.
\item
Assume that $(y^{(n)})_{n\geq 1}$ satisfies (\ref{eqn:sequence}).
Fix $t\in {\bf R}$.
Define 
%
%
\begin{equation}
\label{eqn:eminus}
z^{(n)}\equiv e^{\sqrt{-1}t\log n}\cdot y^{(n)}\quad(n\geq 1).
\end{equation}
\end{enumerate}

%
%
\begin{prob}
\label{prob:last}
Find a sequence which satisfies (\ref{eqn:sequence})
except above examples.
\end{prob}
\ssfr{Acknowledgments}
The author would like to express his sincere thanks to Kentaro Hamachi
for his lecture on Kac-Takesaki operators.

%
%

%
\end{document}